\documentclass[12pt]{article}
\usepackage{amsmath, amsthm, amsfonts, stmaryrd}

\def\FF{\mathbb{F}}
\def\NN{\mathbb{N}}

\def\QQ{\mathbb{Q}}
\def\ZZ{\mathbb{Z}}

\def\Fqbrack{\FF_q \llbracket t \rrbracket}
\def\Fqgen{\FF_q \llbracket t^{\QQ} \rrbracket}
\def\bv{\vec{v}}
\DeclareMathOperator{\Gal}{Gal}

\def\fp{\frac{1}{p}}

\newtheorem{theorem}{Theorem}
\newtheorem{cor}[theorem]{Corollary}
\newtheorem{prop}[theorem]{Proposition}
\newtheorem{lemma}[theorem]{Lemma}

\begin{document}

\title{Algebraic Generalized Power Series and Automata}
\author{Kiran S. Kedlaya \\
Department of Mathematics, Evans Hall \\
University of California, Berkeley \\
Berkeley, CA 94720-3840
}
\date{September 6, 2001}

\maketitle

\begin{abstract}
A theorem of Christol
states that a power series over a finite field is algebraic over
the polynomial ring if and only if its coefficients can be generated by
a finite automaton. Using Christol's result, we prove that the same
assertion holds for generalized power series (whose index sets may be
arbitrary well-ordered sets of nonnegative rationals).
\end{abstract}

\section{Introduction}

Let $\FF_q$ be a finite field of characteristic $p$. As usual, let
$\FF_q[t]$ and $\Fqbrack$ denote the rings of polynomials and
of power series, respectively.
Christol \cite{christol} (see also \cite{ckmr}) proved that an element
$x = \sum_{i=0}^\infty x_i t^i$ of $\Fqbrack$
is algebraic over $\FF_q[t]$ if and only there is a finite automaton
with input alphabet $\{0, \dots, p-1\}$ and output set $\FF_q$ that,
on input the base $p$ representation of a nonnegative integer $i$,
returns $x_i$.
Christol's theorem extends easily to cover power series in a fractional
power of $t$, also known as Puiseux series.

Over a field of characteristic
0, the ring of Puiseux series is algebraically closed (by which we mean
integrally closed in the algebraic
closure of its fraction field), but this is not the case in characteristic $p$.
On the other hand, there is a ring of ``generalized power series'',
first introduced by Hahn \cite{hahn} in 1907, which is algebraically closed.
Specifically, a generalized power series is a sum $\sum_{i \in I} x_i t^i$
where the index set $I$ is a well-ordered subset of the nonnegative rationals.
(An ordered set is \emph{well-ordered} if every nonempty subset has a least
element.) Denote the ring of generalized power series over $\FF_q$
by $\Fqgen$; not only is this ring algebraically closed,
but the algebraic closure of $\Fqbrack$ therein can be identified explicitly
\cite{me}. Thus it is natural to ask whether Christol's criterion admits
a generalization to generalized power series.

The purpose of this paper is to establish an automata-theoretic
criterion for algebraicity over $\FF_q[t]$ of a generalized power series.
Namely, let $x = \sum_{i \in I} x_i t^i$ be a generalized power series over
$\FF_q$. We say $x$ is \emph{automatic} if:
\begin{enumerate}
\item[(a)] The support $I$ of $x$ is contained in $\ZZ[\fp]$;
\item[(b)] there exists a finite
automaton with input alphabet $\{0, \dots, p-1, .\}$ and output set $\FF_q$
that, on input the base $p$ representation of a $p$-ary rational $i$,
returns $x_{i}$.
\end{enumerate}
In terms of this definition, we have the following theorem.
\begin{theorem} \label{thm:main}
Let $x = \sum_{i \in I} x_i t^i$ be a generalized power series over $\FF_q$
whose support is contained in $\ZZ[\fp]$.
Then $x$ is algebraic over $\FF_q[t]$ if and only if it is automatic.
\end{theorem}
\begin{cor}
Let $x = \sum_{i \in I} x_i t^i$ be a generalized power series over $\FF_q$.
Then $x$ is algebraic over $\FF_q[t]$ if and only if for some (any)
$m \in \NN$ such that $mI \subseteq \ZZ[\fp]$, the series
$\sum_{i \in I} x_i t^{mi}$ is automatic.
\end{cor}

\section{Definitions}
Before proceeding, we fix our terminology relating to finite automata.
A \emph{finite automaton} is a tuple $(I, O, S, i, f, T)$ 
in which $I$, $O$, and $S$ are finite sets (the \emph{input alphabet},
\emph{output set} and \emph{state space}, respectively), $i$
is a member of $S$ (the \emph{initial state}), $f$ is a function from
$S$ to $O$ (the \emph{output map}),
and $T$ is a map from $I \times S$ to $S$ (the \emph{transition map}).
Given a string $s = s_1\cdots s_n$
of elements of the input alphabet, define the sequence $i_0, \dots, i_n$
by $i_0 = i$ and $i_{j+1} = T(i_j, s_j)$ for $i=1, \dots, n$;
we call $f(i_n)$ the \emph{output} of the automaton on $s$.

In computer science, the output of a finite automaton is usually
a boolean variable, and states are classified as accepting or rejecting
depending on whether $f$ maps them to the true or false value. Since
our output space is a general finite field, we modify the terminology
slightly. We say a state $i \in S$ is a \emph{zero state} if $f(i) =0$
and a \emph{nonzero state} otherwise; we analogously classify strings
as zero or nonzero.

\section{Proof of the Theorem}

Recall the classification of algebraic generalized power series
given in \cite{me}, here specialized to the coefficient field $\FF_q$.
For $a \in \NN$ and $b,c \geq 0$,
define the set
\[
S_{a, b, c} = \left\{ \frac{1}{a}(n - b_1 p^{-1} - b_2 p^{-2} - \cdots):
n \geq -b, b_i \in \{0, \dots, p-1\}, \sum b_i \leq c \right\}.
\]
\begin{prop}
For $x \in \Fqgen$, $x$ is algebraic over $\Fqbrack$
if and only if the following conditions hold.
\begin{enumerate}
\item[(a)]
The support of $x$ is contained in $S_{a,b,c}$ for some $a,b,c$.
\item[(b)]
For some (any) $a,b,c$ such that the support of $x$ is contained in
$S_{a,b,c}$, there exist $M$ and $N$ such that
every sequence $\{c_n\}_{n=0}^\infty$ of the form
\begin{equation} \label{eq:index}
c_n = x_{(m- b_1 p^{-1} - \cdots - b_{j-1} p^{-j+1} -
p^{-n}(b_j p^{-j} + \cdots))/a } \quad (n \geq 0),
\end{equation}
with $m \geq -b$ and $b_i \in \{0, \dots, p-1\}$ with
$\sum b_i \leq c$, becomes eventually periodic with period length at most
$N$ after at most $M$ terms.
\end{enumerate}  
\end{prop}

The main ingredient in the proof of the theorem, beyond Christol's result,
is the following lemma.
\begin{lemma} \label{lem:bounded}
  Let $x = \sum x_i t^i$ be a generalized power series over $\FF_q$
whose support is bounded and contained in
$\ZZ[\fp]$. Then the following are equivalent:
\begin{enumerate}
\item[(a)] $x$ is algebraic over $\FF_q[t]$;
\item[(b)] $x$ is algebraic over $\Fqbrack$. 
\item[(c)] $x$ is automatic;
\end{enumerate}
\end{lemma}
\begin{proof}
It is clear that (a) implies (b); we now prove that (b) implies (c).
Suppose $x$ is algebraic over $\Fqbrack$; without loss of generality,
we may assume the support of $x$ is contained in $S_{a,b,c}$ for some
$a,b,c$, and that it is also contained in $[0,1)$. 
We now describe an automaton that shows that $x$ is automatic.
Let $S$ be the set of finite strings of elements of $I = \{0, \dots, p-1\}$
having at most $c$ elements not equal to $p-1$, and having no consecutive
string of $M+N!$ elements all equal to $p-1$; clearly $S$ is a finite set.
We use $S$ plus one dummy state $z$
as the state space for a finite automaton as follows. Take
the empty string in $S$ as the initial state. Given a state $s \in S$ and
an element $x \in I$, define the transition map $T(x,S)$ as follows:
\begin{itemize}
\item if $s = z$, put $T(x,s) = z$ independent of $x$;
\item if the concatenation $s+x$ belongs to $S$, put $T(x,s) = s+x$;
\item if the concatenation $s+x$ ends with $N!$ elements equal to $p-1$,
let $t$ be the string obtained from $s+x$ by removing these $N!$ elements,
and set $T(x,s) = t$;
\item otherwise, set $T(x,s) = z$.
\end{itemize}
Define the evaluation map $f(s)$ as follows. If $s = \emptyset$ or $s = z$,
put $f(s) = 0$. Otherwise, if $s = i_1\cdots i_m$, let
\[
f(s) = x_{i_1p^{-1}+\cdots+i_mp^{-m}}.
\]
Then by the criterion for algebraicity given above, this automaton computes
$x_n$ for all $n \in [0,1)$, so $x$ is automatic.

We next show that (c) implies (b).
We first make some reductions to simplify the
structure of $A$.
First, assume without loss of generality that
all strings ending in 0 are zero. For example, this can be arranged by
doubling the state space of $A$, and having transitions always go to one
copy and all other transitions go to the other copy. Next, assume that
every state is reached by at least one string, by deleting all states
that are not reached. Also, collapse all states that can only lead to
zero states into a single terminal state.

Impose a partial order on the state space of $A$
by setting $i \geq j$ if there exists a sequence of transitions from $i$
to $j$, and put $i \sim j$ if $i \geq j$ and $j \geq i$; in the latter
case, we say $i$ and $j$ are \emph{equivalent}. 
In this notation, there is at most one transition from a nonterminal
state $i$ into
the set of states equivalent to $i$. Indeed, if this were not the case,
there would be two sequences $u$ and $v$ that cause $i$ to transition to
itself; if $t$ is a sequence leading to $i$ and $w$ a sequence leading
from $i$ to a nonzero state, and $u$ precedes $v$ in lexicographic order,
then $tvuw, tuvw, tuuvw, \dots$ is an infinite decreasing sequence of nonzero
strings, contradiction. The same argument shows that if $i$ is a state
and $a, b \in \{0, \dots, p-1\}$ are such that $a$ carries $i$ to an equivalent
state while $b$ carries $i$ to an inequivalent, nonterminal state, then
$a > b$.

Let $m$ be the number of equivalence classes of nonterminal states,
and let $n$ be a number divisible by the sizes of all of the equivalence
classes. We prove that the set of nonzero strings of $A$ is contained
in $S_{p^n-1,b,c}$ for some $b$ and $c$; it suffices to exhibit a constant
$d$ such that for all rationals $r$ corresponding to nonzero strings of $A$,
the base $p$ representation of $(p^n-1)r$ has at most $d$ digits not equal to
$p-1$.
Indeed, each nonzero string $s$ can be written as a concatenation
$s_1 \cdots s_l$ for some $l \leq m$, where each $s_i$ except for $s_l$
consists of a sequence of transitions
within an equivalence class followed by a single transition out of the
equivalence class; $s_r$ is similar but has no final
transition out of the class.

To obtain the base $p$ representation of $(p^n-1)r$, where $r$ is the rational
with base $p$ representation $s$, we write $s$ shifted left by $n$ in one
row, write $s$ in a second row, then subtract. We claim that if some $s_i$ 
contains more than $n$ digits, then all but the rightmost $n+1$ digits
in the copy of $s_i$ in the second row will lie above digits of the difference
equal to $p-1$. Indeed, these digits all lie below identical digits
in the first row. Moreover, the $n$-th digit from the right in the
copy of $s_i$ in the second row causes
a transition from a particular state $a$ to an equivalent state, whereas the
digit above it causes a transition from $a$ into an inequivalent state.
As noted earlier, the latter digit is smaller than the former, so the
subtraction causes carry in that position, and in all of the positions to
the left until we run out of digits of $s_i$ in the second row.
This proves the claim;
we conclude from it
that the number of digits of $(p^n-1)r$ not equal to $p-1$ is
at most $m(n+1)$, so the set of rationals corresponding to nonzero strings
is contained in $S_{a,b,c}$ for $a=p^n-1$ and $b,c$ suitably chosen.

Next, let $B$ be an automaton that generates the coefficients of $y
= \sum_i x_i t^{ai}$.
To show that $y$, and hence $x$,
is algebraic over $\Fqbrack$,
it now suffices to show that there exist
$M$ and $N$ such that every sequence of the form (\ref{eq:index})
 becomes periodic
with period at most $N$ after at most $M$ terms. In fact, this holds
with $M$ and $N$ equal to the number of states of $B$: the $n$-th index
in the sequence can be written as $us_nv$, where $n$ is a string of
$n$ $(p-1)$'s, and the state of $us_n$ becomes eventually periodic in $n$
with head and period length no greater than the number of states.
Thus $x$ is algebraic over $\Fqbrack$.

Finally, we prove (b) implies (a). Suppose $x$ satisfies the polynomial
$\sum_{j=0}^n a_j x^j = 0$ for some $a_j \in \Fqbrack$, not all zero.
Write $a_j = \sum_{i=0}^\infty c_{j,i} t^i$. Then for each $l \in [0,1)$,
the condition that the coefficient of $t^l$ in $\sum_{j=0}^n a_j$ is zero
yields a linear recurrence relation of the form
\begin{equation}
\sum_{i=0}^k \sum_{j=0}^n d_{j,i} c_{j,m+i} = 0
\end{equation}
for $m=0,1,\dots$, plus a finite number of additional linear relations
among the $c_{j,i}$. The set of solutions of (1) is finite dimensional
over $\FF_q$ and stable under left shift
(i.e., under replacing $c_{j,i}$ with $c_{j,i+1}$). Thus some two powers of
the left shift act by the same linear map on the set of solutions; this
implies that each solution is eventually periodic. In particular, there exist
$b_j \in \FF_q[t]$, not all zero, such that $\sum_{j=0}^n b_j x^j =0$,
so that $x$ is algebraic over $\FF_q[t]$. This completes the proof.
\end{proof}

We now proceed to the proof of the main theorem.
We first prove that every automatic series is algebraic over $\FF_q[t]$.
If $x$ is automatic, write $x = \sum_{n=0}^\infty y_n t^n$, where
$y_n = \sum_{0 \leq i < 1} x_{n,i} t^{i}$. The set of all of the $y_n$
is finite, since each one is produced by a single automaton starting in
different states. In fact, we can write $x = \sum_{j=0}^k x_j z_j$,
where $x_j \in \Fqbrack$ is automatic and $z_j$ is one of the $y_n$.
By the lemma, each of the $y_n$, being automatic, is algebraic
over $\FF_q[t]$; by the theorem of Christol, each $x_j$ is algebraic
over $\FF_q[t]$. Hence $x$ is algebraic over $\FF_q[t]$.

Before proving the reverse implication, we record a lemma about 
a certain type of semilinear equation. This lemma is closely related
to the Dieudonn\'e-Manin classification of $F$-crystals, but we will not
pursue this relationship further here.
\begin{lemma} \label{lem:vector}
Let $K$ be a field of characteristic $p>0$, let $q$ be a power of $p$,
and let $A$ be an $n \times n$
matrix with entries in $K$. Let $\bv$ be an $n$-vector with entries
in some field containing $K$ such that $A\bv^q = \bv$. Then the entries
of $\bv$ are algebraic over $K$.
\end{lemma}
\begin{proof}
First note that the set of vectors $\bv$ with entries in any
field $L$ containing $K$ such that $A \bv^q = \bv$ is an $\FF_q$-vector
space of dimension at most $n$. Indeed, let $\bv_1, \dots, \bv_m$ be
a maximal $K$-linearly independent set of solutions, and let $\bv
= c_1 \bv_1 + \cdots + c_m \bv_m$ be any other solution;
then $c_i^q = c_i$ for $i=1, \dots, m$, so $\bv$ is actually in the
$\FF_q$-span of $\bv_1, \dots, \bv_m$.

In particular, the number of solutions of $A \bv^q = \bv$ is finite
in any field containing $K$. If the entries of $\bv$ were not all algebraic
over $K$, there would be an extension $L/K$ such that the orbit of $\bv$
under $\Gal(L/K)$ is infinite. Thus the entries of $\bv$ must be algebraic.
\end{proof}

Now suppose $x \in \Fqgen$ is algebraic and has support in $\ZZ[\fp]$.
Write $x = \sum_{n=0}^\infty t^n
y_n$, where each $y_n$ is supported on $[0,1)$. By Lemma~\ref{lem:bounded},
the set of all of the $y_n$ is finite, and each $y_n$ is algebraic over
$\FF_q[t]$ and automatic. Let $z_1, \dots, z_r$ be a basis
of the $\FF_q(t)$-span
of $y_i^{q^j}$ for $i, j \geq 0$ consisting of series with bounded support,
and write $x = \sum_{i=1}^r c_i z_i$ with $c_i \in \Fqbrack$.

Since $x$ is algebraic over $\ZZ[\fp]$, there must be a relation of the
form $\sum_{i=0}^m a_i x^{q^i} = 0$ for some $a_i \in \FF_q[t]$ not all 
zero. There is no loss of generality in assuming $a_0 \neq 0$ (by replacing
$x$ with a power thereof). Substituting $x = \sum_{i=1}^r c_i z_i$
into the relation and rewriting each $z_i^{q^j}$ as an $\FF_q(t)$-linear
combination of the $z_i$, we get a linear relation of the form
$\sum_{i=1}^r f_i(c_1, \dots, c_r) z_i = 0$. More precisely,
there exists $N>0$ such that each $f_i$ can be written
as $c_i + \sum_{l=1}^r \sum_{j=1}^N d_{i,l,j} c_l^{q^j}$ for some
$d_{i,l,j} \in \FF_q(t)$.
Now consider the system of equations
\begin{align*}
  e_{i,j} &= e_{i,j+1}^q \qquad (i=1, \dots, r; j=1, \dots, N-1) \\
  e_{i,N} &= \sum_{l=1}^r \sum_{j=1}^{N} d_{i,l,j} e_{l,N+1-i}^q \qquad (i=1,
\dots, r).
\end{align*}
On one hand, this system has the solution $e_{i,j} = c_l^{q^{N-j}}$.
On the other hand, it has the form $A \bv^q = \bv$ for a suitable matrix
$A$ over $\FF_q(t)$, so by Lemma~\ref{lem:vector} any solution has entries
algebraic over $\FF_q(t)$. We conclude that the $c_i$ are algebraic over
$\FF_q(t)$.

Now by the theorem of Christol, each of the $c_i$ is automatic.
Since $z_i$ has bounded support and is automatic, and $c_i \in 
\Fqbrack$, $c_i z_i$ is also automatic. Finally, we conclude
$x = \sum_i c_i z_i$ is automatic. This completes the proof of 
Theorem~\ref{thm:main}.

\subsection*{Acknowledgments}
Thanks to Bjorn Poonen for bringing the work of Christol to the author's
attention. The author was supported by an NSF Postdoctoral Fellowship.

\end{document}